# Turing-Church thesis, constructve mathematics and intuitionist logic

Antonino Drago
University "Federico II" of Naples – drago@unina.it

**ABSTRACT**   At a first glance the Theory of computation relies on potential infinity and an organization aimed at solving a problem. Under such aspect it is like Mendeleev's theory of chemistry. Also its theoretical development reiterates that of this scientific theory: it makes use of doubly negated propositions and its reasoning proceeds through *ad absurdum* proofs; a final, universal predicate of equivalence (of all definitions of a computations) is translated into an equality one, and at the same time intuitionist logic into classical logic. Yet, the last step of this development of current theory includes both a misleading notion of "thesis" and intuitive notions (e.g. the partial computable function, as stressed by some scholars). A program for a rational re-construction of the theory according to the theoretical development of the above mentioned theories is sketchy suggested.
**Keywords**: Constructive mathematics, Intuitionist logic, Turing-Church thesis, Theory of computation, Doubly negated propositions, *Ad absurdum* proofs, Translation of the kind of logic, Partial computable function, Rational re-construction of the theory

## 1. Turing-Church's thesis, constructive mathematics and intuitionist logic

The question of which is the relationship between through Turing-Church thesis (TCT) and constructive mathematics leads to investigate on one hand, the theory computation (TC). Unfortunately, there is no common agreement whether TC is a scientific theory; also because ultimately there is no common agreement about the definition of the notion of a scientific theory (Stanford 2017). On the other hand, one has to investigate on constructive mathematics whose supporters declare that it is classical mathematics but according to intuitionist logic.(Bridges 2018) Hence, the question concerns also the relationship between TCT and intuitionist logic (IL).

In the past the question of the meaning of TCT has been examined under classical logic (CL). I attribute the inconclusiveness of this debate to this constraint, because, as it will be manifest in the following, TCT belongs to IL.

## 2. The theoretical development of computer theory: A basic problem

First, let us consider TC as a theory. It is comparable with Mendeleev's theory of chemistry. Mendeleev's textbook on Chemistry starts by posing two basic problems. The first one is of a "qualitative nature": "What is an element?" (Mendeleev 1906, pp. 20-21). In TC this kind of problem is easily recognised as the following one: "What is a computation?". In order to solve this problem, Mendeleev's theory first solves a subordinate problem of a "quantitative nature": "How much chemical reactions one has to perform for recognizing an element?". In TC we recognize the quantitative problem as the following one: "How much kinds of computation notions one has to discover in order to define a computation?" In Chemistry, after having examined the results of a number of chemical reactions and moreover having attributed a valence number to each element (that is the relational import of the theory), Mendeleev established the periodic table of elements which represents the unitarian nature of all possible elements of matter. In such a way he answered the qualitative problem: an element is a component of his table.(Drago 2014) After having examined several formal notions of computability and moreover proved that they all are mutually equivalent (that is, the relational import of the theory), TC, through Turing-Church's thesis (TCT), establishes the equality of whatsoever kind of formal computation with the intuitive notion of computability. That answers the qualitative problem of current TC.

## 3. The theoretical development of theory of computation: Doubly negated propositions

Moreover, both theories make use of propositions of a particular kind, i.e. doubly negated propositions which are not equivalent to the corresponding affirmative propositions since the semantic contents of the latter ones lack of evidence within the real world (DNPs). In last century the scholars of mathematical logic achieved a crucial result; i.e. the validity or not of the double negation law represents the best discriminating mark between classical logic and most non-classical logics, above all intuitionist logic (Prawitz and Melmnaas 1968; Grize 1970, pp. 206-210; Dummett 1977, pp. l7-26; Troelstra and van Dalen 1988, pp. 56ff.). Hence, this failure of the double negation law in the cases of the DNPs qualifies all them as belonging to non-classical logic, in particular, intuitionist logic.[1]

TC's includes a modal[2] adjective, 'partial' in the same subject of study: "partial computable functions". Notice that a modal word is equivalent to a DNP; e.g. possible = it is not the case that it is not". Also in formal terms, modal logic may be translated by means of its S4 model into intuitionist logic. (Hughes and Cresswell 1996, pp. 224ff.) Actually, the meaning of this word "partial" is more clear when it is replaced by an equivalent DNP which deals with only what we can objectively know; hence, not its unverifiable undefined values, but only the determinations of a final results of computing a function; hence, the correct translation of "partial" in objective terms is the following one: "not never defined".

In addition, instances of DNPs occur in an important TC's point of its development. Although ignoring the logical features of a DNP, several authors have presented TCT through propositions of this kind. For instance, Stephen Kleene: "… it cannot conflict..." (Kleene 1952, pp. 318-3l9). Alan Turing makes use of two modal words, "could" and 'naturally': "the [Turing machine's] computable numbers *include* all numbers which could naturally be regarded as computable";(Turing 1936, p. 249) Alonso Church avoids a DNP by means of a turn of words: "it is thought to correspond satisfactorily..." (Davis 1965, p. 90)[3] the same is done by Martin Davis: "we have reason to believe that...",(Davis et al. 1995, pp. 68-69; without specifying this reason); moreover, its word "believe" is a subjective word representing a modality of truth; it actually means the DNP: "it cannot be false that it is…" Kurt Goedel has called TCT a "heuristic principle". (Davis 1965, p. 44) Equivalently, Emily Post: "A working hypothesis" (Davis 1965, p. 291).

As a consequence of the use of all these DNPs, TC makes an essential use of non-classical logic. This is the real reason for calling the subject of this theory by means of a modal word, 'computability' (see the title of Davis et al. 1995); it means that the theory deals with what *may* be computed.

Notice that intuitionist logic cannot govern a deductive theory since each DNP does not oppose to the corresponding affirmative proposition (to the negative one either); hence, it does not give a true/false mirror opposition that assures certainty to each step of a logical derivation. As a consequence, classical logic governs AO theories, intuitionist logic governs a theory as the chemistry, which is based on the search for a new method solving a given problem.[4] Notice that also some other theories share the previous features: L. Carnot's mechanics (1783), Lagrange' mechanics (1878), Sadi Carnot's thermodynamics (1824), Galois' theory (1832), Klein projective theory of geometries (1881), Einstein's first paper on quanta (1905), Kolmogorov's foundation of intuitionist logic (1932), Markov foundation of computable numbers (1962).(Drago 2012)

I call the organization of all these theories a problem-based one.

## 4. The theoretical development of Theory of computation: *ad absurdum* arguments

---

[1] Independently, Grzegorczyk (1964) proved more in general that the production of new results by experimental science may be formalized through propositions belonging to intuitionist logic, that is, a logic making use of DNPs. It is the theoretical physicist who truncates the approximations values of experimental results in order to obtain idealistically exact results; ie he translates the rational numbers of measurements into real numbers and hence intuitionist logic into classical logic.

[2] A modal word will be underlined with dots.

[3] Most textbooks make use of the word "equivalence", which however is not independently defined by them. Notice that the same word is used also by most Thermodynamics' textbooks illustrating the first principle, concerning all kinds of energy. Actually, this word means: "It is not true that it is non equal...", i.e. it is a DNP.

[4] A comparative analysis of a lot of scientific theories formally proves the last propositions. (Drago 2012)

An author of a PO theory links through IL a DNP to another DNP in order to constitute an *ad absurdum* argument (AAA) (revealed within the original text of the theory by words like e.g. "otherwise… not…"), whose conclusion is no other than a DNP. Moreover, a DNP concluding an AAA may work as premise for a next AAA, so that such arguments constitute a logical chain.

Also TC presents a great number of *ad absurdum* arguments, as a glance at TC's textbooks shows; on the contrary, the usual mathematical theories, being based on classical logic, seldom make use of them.

In general the conclusion of this chain of arguments is a universal predicate; it represents a guess for solving the given basic problem and all correlated problems of the theory. As a matter of fact, the predicates concluding the above mentioned theories all are equivalent to the main thesis A of the square of opposition (see the table in Dummett 1977, p. 29). At this point the theoretical development a PO theory proceeds in an alternative way to the previous use of only logical arguments; it leads to make use of mathematical formulas which, by unavoidably including an exact equality symbol, belong to classical logic. This goal is achieved by translating the doubly negated, conclusive predicate into the corresponding affirmative predicate, which only can represent a mathematical formula; that amounts to both a translation of the kind of logic, from CL to IL, and a change of theory's organisation, from previous PO to AO, which is the suitable organization for linking together mathematical derivations based on exact equalities.

The above translation of the universal doubly negated predicate may be justified in the following way: since at this moment the predicate is supported by the previous formal reasoning including some AAPs, the author feels himself justified in translating it into its corresponding affirmative predicate.

Notice that this change results from an (implicit) application of Leibniz' principle of sufficient reason (= PSR) – whose antecedent is itself a universal doubly negated predicate ("Nothing is without reason"[5]), whereas the consequent is an affirmative one ("There exists a reason") which belongs to classical logic. This translation changes the intuitionist version of thesis A of the intuitionist square of opposition (see Dummett's table) into the classical thesis A. It is easy to recognize through the same table that the logical translation of the intuitionist thesis A into the classical one amounts to translating the entire kind of predicate logic, from the intuitionist one to the classical one. Under this light, the application of the PSR appears as the inverse translation of the well-known translation through double negations from classical logic into the intuitionist one (Troesltra van Dalen 1988, pp. 56ff.).[6]

Is there a conclusive doubly negated predicate within TC? After having proved that the formal definitions of computation all are mutually equivalent - in the sense they all define the same set of functions (Odifreddi 1986, pp. 87-101) – TC's textbooks compare these definitions of computation, $C_f$, with the intuitive notion, $C_i$.

This final step of TC's theoretical development appears similar to that of a PO theory. Yet, let us more accurately inspect this step. As a matter of fact, TC's textbooks present the relationship between the two notions as a possible "equivalence" of $C_f$ and $C_i$ ; which means the following universal predicate: "It is not true that the intuitive notion of a computation is not the formal notion of a computation", $\neg(C_f \neq C_i)$; i.e. a DNP. Their aim is to state the corresponding affirmative predicate one, ie the equality, $C_f = C_i$; which is the starting point of a subsequent deductive development of the theory. As a fact, this equality is stated by appealing to a so-called "thesis", precisely TCT.

---

[5] "Or - Leibniz continues - everything has a cause,… although we are not always able to show it." (Leibniz 1714) Leibniz' last proposition admirably stresses that the former DNP cannot be formally equated to its affirmative proposition, exactly what modern non-classical logic teaches.(Drago 2017)

[6] All that agrees with what Leibniz maintained: there exist two basic logico-philosophical principles, the principle of non contradiction and the principle of sufficient reason. From all in the above we may conclude that the former principle is associated to the AO of a scientific theory, whereas the latter one to the PO.

This word is usually intended as the declaration of a subsequent theorem; which however TC's theorists cannot obtain for the plain reason that from an intuitive notion one cannot formally derive the corresponding formal notion. We have to conclude that the proof of this "thesis" is a strange "theorem", which, although seemingly plausible, cannot exist. Being prevented to prove this "equality", textbooks, make use - as we saw in previous sect. -, of turns of words inducing readers to think that the thesis is in some way formally valid.

As a historical fact, the belief in this "thesis" fortunately produced an exciting achievement – the convergence of all formal notions of computability into one.[7] After this achievement, all its possible premises seemed as justified or seemed in course to be justified by some imminent achievements of the research on foundations. But at present textbooks cannot still exhibit a certain proof of TCT.

In conclusion of previous sect.s 2-4, one may recognize some formal steps of the theoretical development of a PO theory inside the theoretical framework of TC; but the last step of TC's development is odd with respect to that of each PO theory, so much to make dubious the characterization of TC as a PO theory.

### 5. Intuitive notions and formal notions

TCT includes a notion– ie "thesis" - whose meaning is ambiguous because, although it alone represents an intuitive notion, is used by TC as pertaining to a formal argumentation. TCT is a strange element of present TC, which as a whole is instead a highly formalized theory. To some scholars the word "thesis" of TCT appears inappropriate for denominating a scientific result.[8] Following Galton's words lucidly characterise the situation:

> We thus seem to have the curious circumstance of a corpus of exact mathematical results being held up as evidence for a thesis which is impossible to state clearly while maintaining its distinctness from those results.(Galton 1996, p. 141)

The inclusion by TC of this odd word "thesis" - ie a notion which being no more than intuitive in nature, yet plays the role of "proving" the formal equality $C_f = C_i$ - leads to ask the following question: What means to built a scientific theory in order to define an intuitive notion? A plain answer is the following one: if at a certain step of its development a theory formally obtains the wanted result, from this time on the notion at issue is no longer an intuitive one, because it was translated through formal means into a formal element of a scientific context;[9] alternatively, if the theory cannot define it as a formal notion, the trivial conclusion is that the theory is merely unsuccessful, exactly what we have to conclude about TC, which manifestly cannot prove TCT.

If this notion has to be preserved by the theory as an intuitive notion, it represents a stumblock to a scholar inspecting TC's foundations. Surely, in the case of a theory including, as a basic proposition of the theory, a "unproved thesis" we cannot organize it as an AO, since by starting from an intuitive, hence inaccurately defined notion, no accurate logical inference of an axiomatic theory is possible.

Otherwise, may be TCT considered as a program of research for solving the following problem: is $C_f = C_i$ true? This suggestion actually addresses to consider TC as a PO theory aimed at solving this problem. According to the model of a PO theory, not before having applied the PSR for translating the conclusive doubly negated predicate into the corresponding affirmative predicate, TCT obtains an affirmative predicate from which all the consequences may be deduced according to classical logic. Thus, TC can translate the *equivalence* (= <u>not</u> <u>difference</u>) $\neg(C_f \neq C_i)$ of the intuitive notion of computability and the formalized ones, into an *equality* $C_f = C_i$ not by applying an intuitive notion, as the "thesis" (ie TCT) is, but the PSR.

---

[7] So much that the usual TCT was exalted by Goedel as leading to define 'computability' as an "absolute notion" (Goedel 1990, p. 150). In my opinion this word 'absolute' inaugurated a new scientific idealization without any clear explanation.

[8] Present TC seems an application of a Suppes' philosophical attitude which accepts that a rigorous axiomatic of a physical theory cannot exist, because an axiomatic should include also connection axioms with reality, which however are impossible, given the informal nature of the reality. "A Suppes' axiomatic" is constituted by no more than a mathematical predicate applied to a *mathematical model* of reality.

[9] This is the case of Information Theory. There, a merely mathematical definition of "information" is cursorily justified and then from it all consequences are drawn, by disregarding any question of the adequacy of the mathematical notion of "information" with the corresponding intuitive notion, since the latter notion, by apparently including many semantic meanings, has surely much more wide import than that of the former one.

Yet, an application of PSR to a predicate is correct when PSR is applied to a decidable predicate resulting from an AAA.(Drago 2012; Drago 2017) But an AAA concluding TC's final predicate, which includes an intuitive notion $C_i$, can refer to no more than an informal rationality (e.g.: "Otherwise, our mind errs"). In more precise terms, previous predicate is not decidable. Hence, no correct translation of this predicate of equivalence into an equality predicate is possible through PSR.

In conclusion, present TC, by including the intuitive notion of computable function, can be organized according to neither the model of an AO, nor the model of a PO. As a matter of fact, along many decades scholars unsuccessfully tried to obtain a specific organization of TC.[10]

### 6. A crucial choice: the inclusion of actual infinity into TC

How it was possible that TC has entered into such an apparently blind alley? The origin of this problem is apparently constituted by the intuitive notion, $C_i$ which plays an essential role within this scientific theory.

Let us recall that Koyré's celebrated analysis of Galilei's work underlined that the Italian scientist shared a Platonist conception of Mathematics, according to which Mathematics represents a language written within the reality (Koyré 1978): "The great book of the Universe is written in squares, circles and triangles,…". (Galilei 1623, pp. 121-122) As a consequence, the word 'intuition' means that human mind enjoys an ontological insight within the reality. Instead, a century and half before Galilei, Cusanus had suggested that Mathematics is a construction of our mind; we cherish it because "nothing is more certain than our Mathematics", being this activity rigorously founded on the principle of non-contradiction.(McTighe 1970). This conception leads to assume the opposite philosophy of Galilei's relationship mathematics-physics as an ontological intuition at reality because human intuition may be true or wrong; it is true when an intuitive notion is generated by either an (at least partially) formal notion or experimental facts; in other words, it is true when it subjectively summarizes a formal notion in a convenient way. As a matter of fact, each advancement in the history of Mathematics has suggested an intuitive notion of what in a previous time was considered a fantastic or even an absurd notion.[11] More in general, by following the mathematical progress human intuition is widening ever more its realm, provided that its relies on well established formal notions. In short, in science intuition comes after the formal notions.

From this viewpoint a completed scientific theory has to reject all intuitive notions; in particular, TC has to intend the word 'computation' in objective terms, not as the intuitive notion of "computability". Instead, present TC, being committed to an intuitive notion $C_i$ as an a priori formally undefined notion, assumes a Platonist view of Mathematics, truly the dominant view in present philosophy of science.[12] As a consequence, although presenting an apparently strictly operative theory, TC's textbooks include at least this philosophical premise.

---

[10] In order to clarify the theoretical situation of present TC let us again compare it with the theoretical developments of a manifest PI theory rejecting any Platonist notion, Mendeleev's chemistry. Modern chemists oriented their theory by rejecting alchemists' belief that an intuitive element, the *hyle*, could be transformed into each of the several elements of matter and *viceversa*. Two centuries after, TC's theorists based their theory on the intuitive notion of computability which was declared equivalent to the formal notions of computability in mathematical and physical (Turing machine) terms. Thus, within TC the intuitive notion $C_i$ plays the same theoretical role played by the notion of *hyle* within alchemy, a theory which at present is dismissed because it is not an entirely experimental theory and the notion of *hyle* is falsified by the lack of a transmutation between eg Lead into Gold.

[11] For instance, a millennia long experience of Euclid's geometry generated within human mind a "geometrical intuition" of e.g. the single point shared by two crossing curves, the single point shared by a curve and its tangent, etc. Projective geometry generated an intuitive notion of the point at infinity, that previous painters had badly and obscurely perceived. And so on.

Moreover, TC assumed a Platonist attitude also about the intuitive notion of infinity, which formally has to be distinguished in two formal notions, i.e. the potential infinity (PI), on which relies constructive mathematics, and the actual infinity, on which relies classical mathematics.

In past times, the dominant group of mathematicians held a prejudice with respect to TC; being this theory very useful for practical scopes and being its kind of mathematics very near to the more simple mathematical system, arithmetic, its theoretical import was considered of little relevance.[13] In order to oppose this prejudice, TC's scholars wanted to enlarge as most as possible the theoretical import of their theory, in the aim at achieving as most as possible theoretical relevance.

History of TC seems a rushing to ever more impressive results. In particular, let us sketchily consider the history of recursive functions. As it is well known, elementary recursive theory, relying on the notion of PI, originated through the celebrated Goedel's paper of the year 1931. Then Ackerman invented a recursive function which is different from all elementary recursive ones; it essentially includes AI. After this result, TC's theorists developed without scruples a general recursion theory relying on AI as well as all the idealistic notions which it suggests.

In particular, they also showed that the set of the general recursive functions, relying on AI, is the same of the set of partial computable functions. The word 'partial' means that we do not know for which numbers a function of such a kind is defined.[14] Sundholm has stressed that this notion is inappropriate to TC

> However, strictly speaking, *partial function* is an oxymoron. The adjective `partial' acts as a [structural] *modification* [of its subject, ie the noun function] that takes us out of [the usual definition of] functions, rather than as a qualifying property among functions. A 'partial function' is [actually] no function, since it is not defined for every element of the domain. The syntactic form "partial function' is misleading. Instead, in recursion theory, one could better speak about recursively enumerable functional relations, whether total or not. The reading 'function that is partial recursive', on the other hand, would appear to indicate a (total) function that for some reason is not fully recursive. but only partially so. (Sundholm 2014, p. 3)

One may provocatively ask: Being usual scientific theories based on facts, already occurred and always repeatable, why TC is not a theory of the already computed functions? Is this theory interested in knowing future events - whether a function will obtain a result or not -, instead of present time results? By investigating about partial functions, ie functions which may obtain results, is TC perhaps interested in modalities, instead of hard facts? I did not find out answers within textbooks. The only certainty is that the introduction of AI into TC leads to allow an incontrollable notion of function and as a consequence to obscure the subject of an otherwise operative theory.

---

[12] Popper has introduced a criterion of falsificability for deciding whether a theory is not a scientific one. No theory that includes an intuitive notion satisfies this criterion of falsificability. Eg, in Thermodynamics the notion of total energy of a system, $U$, is not falsifiable by experimental means because we lack of measurements apparatuses for the kinds of energies we will discover in the future; but the balance of this magnitude, $\Delta U$, which concerns presently known kinds of energy, is falsifiable; the theory makes use of only $\Delta U$. Also TC's intuitive notion of equality in the formula $C_f = C_i$ is not falsifiable in the direction $C_f \to C_i$, because for each formal definition of computable function one may creatively invent a subjective corresponding notion (instead the inverse direction $C_i \to C_f$ may be falsifiable because one can invent an intuitive notion of computable function which is not represented by the formal ones, as the historical case of Ackerman's invention of a non primitive recursive function proves). In conclusion, the presence of intuitive notions within a scientific theory makes it a unfalsifiable theory, ie, in Popper's language, a metaphysical theory.

[13] This prejudice compares well with that shared by the mainstream of theoretical physicists with respect to Thermodynamics, called by them a merely "phenomenological" theory.

[14] TC's generalization of the usual notion of a completely defined function into a partial one apparently parallels the generalization that a theoretical physicist performs on each result he receives from an experimental physicist. While this experimental result is a truncated number - hence a rational number with a finite number of digits, hence it is at all known -, the former physicist, in order to make easier through calculus his calculations on it, extrapolates it into a real number; which, being composed by an infinite number of non-periodical digits, as say π, is not completely known (if not in an idealistic, Platonist way); in such a way the dense set of rational numbers is changed into the compact set of real numbers. Yet, whereas a theoretical physicist verifies his calculations performed through real numbers, by comparing them with the rational numbers of measurement data (within an acceptable range of approximation), TC's theorist does not know whether his generalization into partial functions is confirmed or not by hard facts (surely not by the mutual equality of all notions of a computable function - the intuitive one too -, an equality which actually is supported by a claimed "thesis" which lack of a formal proof).

In order to acquire more light on this subject let us notice that the attitude of TC's theorists of considering as most as possible functions, even the partial ones, seems repeat the attitude of 19$^{th}$ Century mathematicians of calculus, wanting to embrace all possible real functions - also pathological ones - under an intuitive notion of a "law", or a "correspondence", between two intuitively conceived "sets". Their original aim was to find out through these new functions the most general definitions of the two basic operations of calculus, ie derivative and integral. After two centuries of an enormous inventive effort, someone correctly remarked that the variety of pathological real functions (even if bounded to be monotone functions! Zamfirescu 1981) is so great that they make impossible any universal definition of both derivative and integral; so that at present time the pathological functions are to be called "monsters", to be removed from the corpus of Mathematics of a working mathematician.(Feferman 2000)

After half a century from their introduction, have we to do the same move with the partial computable functions? A scholar (Turner 2004) uggested this move since he contested TC's scholars decision of basing the theory upon partial computable functions. He proved that this decision is not unavoidable. By choosing only total functions he suggested an innovative programming language ("Total functional programming").[15] This result motivated him to write the following plea for making use of total functions only:

> There is a dichotomy in language design, because of the halting problem. For our programming discipline we are forced to choose between
> A) Security - a language in which all programs are known to terminate.
> B) Universality - a language in which we can write
>    (i) all terminating programs
>    (ii) silly programs which fail to terminate
>    and, given an arbitrary program we cannot in general say if it is (i) or (ii).
> Five decades ago, at the beginning of electronic computing, we chose (B). If it is the case, as seems likely, that we can have languages of type (A) which accommodate all the programs we need to write, bar a few special situations, it may be time to reconsider this decision. (Turner 2004, p. 767)

According to this suggestion, current TCT, concerning the set of partial computable functions, is too powerful since it concerns a too large set of functions with respect to the operatively manageable ones. (Turner 2006) In conclusion, it is realistic and (also more appropriate) to purge TC of partial computable functions, introduced as a consequence of the extension of primitive recursive functions to the general ones (relying on AI); that means to bound TC to make use of PI only.

**7. Towards a rational re-construction of TC according to both the model of a PO theory and constructive mathematics**

A possible re-formulation of TC may no longer include the words "computable" and "computability", but only the word "computed", exactly as whatsoever experimental scientific theory refers to comprovable facts; in other terms, TC has to refer to not a subjective, modal notion, but only objectively defined notions, i.e. total functions.

Let us start this re-construction of TC as a PO. The first step of the theoretical development of a PO theory is the declaration of its basic problem. This problem is not the ontological one: "What is a computation?" because a computation is not a material or idealistic thing, but a process; instead it is the quantitative question "How many notions of computation there exist?"; or better, it is the henological question "How unify the multitude of formal definitions of computation?"

---

[15] Notice that the author follows Goedel's *Dialectica* interpretation, which includes idealistic elements.

In addition, let us re-construct TC through constructive tools only. Hence the theory includes only the elementary recursive functions, without dealing with partial computable functions.

About the subsequent development, what has been illustrated in Sect.s 2-4 showed that the theory of TC can be organized as PO theory, apart the last step of its development, ie the application of PSR. It Is possible to recover this step? Does the wanted re-formulation include an analog of the usual TCT? Yes, it is again a 'thesis', yet with the usual meaning, i.e. it is the declaration of a subsequent theorem; which has to prove, through constructive mathematical tools, the equivalence between any pair of formal notions of *mathematical* computation (elementary recursion theory, constructive theory of Diophantine equations, appropriately bounded Lambda Calculus, etc.). I do not know whether such equivalence proofs exist; but at a first glance they seem possible. In addition, I have no reason for excluding that they make use of AAAs which substantiate the characteristic reasoning of a PO theory. Rather, they seem unavoidable because each times are compared to at all different formalisms which cannot be deduced through constructive tools. one from another

The re-formulation continues by tackling the equivalence between on one hand all mathematical notions of computation $C_M$ and on the other hand the physical notion of a Turing machine computation, TC. A priori the latter notion may be essentially different and even incommensurable with the previous one, since whereas $C_M$ is a notion built by our mind, TC includes in an essential way the relationship with reality; between mind and reality one can only guess an equivalence; of course, stated through a DNP. In current TC this equivalence is already called "the physical TCT", but now it is bounded to comparing total Turing functions and total computation functions; moreover it is clearly enounced as a DNP: "It is <u>impossible</u> that Turing computation gives different results from the mathematical computations"; or, in other temrs: "The two kinds of computation are equivalent"..

At this point one can correctly apply the PSR to this DNP, since the two requirements of PSR's applicability are satisfied: the predicate is clearly decidable. Moreover it comes from an AAA; because this one addresses the question if a Turing machine - contrarily to a mathematical process of computation, whose nature of a construction by human mind is a priori based on the principle of non contradiction -, can generate through computations a contradiction. One easily proves through an AAA that surely it cannot; otherwise the physical reality would represent a contradiction within the millennial calculations by mankind (this principle is similar to the methodological principle of the impossibility of perpetual motion in theoretical physics).

Then the application of PSR to the physical TCT obtains the corresponding affirmative proposition: "Both physical and mathematical notions of computation define the same set of functions". From it the next part of the PO theory may be drawn in a deductive way.

In particular, the intuitive notion of "computability" no longer pertains to the re-formulated theory; this notion has to merely conform to the two formal notions that in the above have been stated to be equivalent.

In retrospect, we see that the usual TC actually generalized the new "physical TCC" in a Platonist way by including the intuitive notion of computable function too and moreover by referring to idealistic notions (a "thesis" without a proof, proofs through mathematical techniques of non-constructive, idealistic nature, because they rely on actual infinity, etc.). Let us remark that the rejection of intuitive notions – as computability, intuitive notion of computation, etc. -, does not dismiss real contents. Rather, we have gained that in this way the intuitive notions of physical and mathematical computation are supported by a scientific, constructive theoretical scheme.

In conclusion, the introduction of the IL led to deny any theoretical role to intuitive, subjective notions, which before have been introduced for absorbing TC's properties belonging to IL into intuitive, vague notions.

**Bibliography**

Bridges D.S. (2018), « Constructive mathematics », in N.E. Zalta (ed.), *Stanford Encyclopedia of Philosophy*, https://plato.stanford.edu/entries/mathematics-constructive/

Davis M. (1965), *The Undecidable*, New York: Raven.

Davis M. et al. (1995), *Computability, Complexity, and Languages. Fundamentals of Theoretical Computer Science*, New York: Academic Press.

Drago, A. (2012), "Pluralism in Logic: The Square of Opposition, Leibniz' Principle of Sufficient Reason and Markov's principle", in *Around and Beyond the Square of Opposition*, Béziau, J. –Y. and Jacquette, D. (eds.), Basel: Birkhaueser, 175-189.

Drago A. (2014), "Il ruolo del sistema periodico degli elementi nel caratterizzare la chimica classica come teoria scientifica", *Epistemologia*, 37 (2014) 37-57.

Drago A. (2017), "A Scientific Re-assessment of Leibniz's principle of Sufficient Reason", in R. Pisano et al. (eds.) *The Dialogue between Sciences, Philosophy and Engineering. New Historical and Epistemological Insights. Homage to Gottfried W. Leibniz 1646-1716*, London: College Publications, pp. 121-140.

Dummett M. (1977), *Elements of Intuitionism*, Oxford: Clarendon Press.

Feferman S. (2000), "Mathematical Intuition vs. Mathematical Monsters", *Synthese*, 125(3), pp. 317-332.

Galilei G. (1623), *Il Saggiatore*, Roma: G. Mascardi.
<https://play.google.com/store/books/details?id=-U0ZAAPAYAAJ>.

Galton A. (1996),"The Church-Turing thesis: Its nature and status", in Millikan P.J.R., Clark A., *Machines and Thought*, Oxford: Claredon, 137-164.

Gödel K, (1990), *Collected Works*, vol. II, Oxford: Oxford U.P..

Grize J.-B. (1970), « Logique », in *Logique et la connaissance scientifique*, in *Encyclopédie de la Pléyade*, Piaget J. (ed.), Paris: Gallimard, 135-288.

Grzegorczyk A. (1964), « Philosophical plausible formal interpretation of intuitionist logic", *Indagationes Mathematicae*, 26, 596-601.

Hilbert D. (1967, orig. 1926), "On the infinity", in *From Frege to Goedel,* van Heijenoort J. (ed.), Harvard: Harvard U.P., 376.

Hughes G.E. and Cresswell M.J. (1996) *A New Introduction to Modal Logic*, London: Routledge.

Jona M. (1984), "What is Energy?", *Physics Teacher*, **22,** 6.

Koyré A. (1978), *Galileo Studies*, Hassocks: Harvester Press

Kleene S. C. (1952), *Introduction to Metamathematics,* Princeton: Van Nostrand.

Lavoisier A. L. (1862-92), *Oeuvres de Lavoisier*, Paris: Imprimerie Imperiale, t. 1.

Leibniz G.W. (1710), *Theodicy*, New York: Cosimo, 2009.

Markov A.A. (1962), "On Constructive Mathematics," *Trudy Math. Inst. Steklov*, **67,** 8-14; Engl. tr. *Am. Math. Soc. Translations,* (1971) **98** (2) 1-9.

McTighe T. (1970), "Nicolas of Cusa 's philosophy of science and its metaphysical background", in Santinello G.(ed.), *Nicolò Cusano agli inizi del mondo moderno*, Firenze: Sansoni, pp. 317-338.-

Mendeleev D. I. (1906), *Principles of Chemistry*, London: Longmans.

Odifreddi P. (1989), *Classical Recursion Theory*, Amsterdam: North-Holland.

Prawitz D. and Melmnaas P.-E. (1968), "A survey of some connections between classical intuitionistic and minimal logic", in *Contributions to Mathematical Logic,* Schmidt H. A., Schütte K. and Thiele H.-J. (eds.), Amsterdam: North-Holland, 215-229.

Stanford K. (2017), "Underdetermination of Scientific Theory", in N.E. Zalta (ed.), *Stanford Encyclopedia of Philsoophy*, https://plato.stanford.edu/entries/scientific-underdetermination/

Sundholm G. (2014), "Constructive recursive functions, Church's thesis and the theory of Brouwer creating subject.", in Dubucs J., Bourdau M., *Constructivity and Computability in Historical and Philosophical Perspective*. Berlin: Springer, 1-35.

Troelstra A. and van Dalen D. (1988), *Constructivism in Mathematics*, Amsterdam: North-Holland, vol. 1.

Turing A. (1936), "On computable numbers, with an application to the entscheidungs-problems", *Proc. London Mathematical Society*, Ser. 2, 42, pp. 230-265.

Turner D. (2004), "Total function programming", *J. of Universal Computer Science*, vol. 10, no. 7, 751-757.

Turner D. (2006), "Church's thesis and functional programming", in Olszewski A. et al. (eds.), *Church's Thesis after 70 Years,* Frankfurt: Ontos.

Zamfirescu F.S. (1981), "Most monotone functions are singular". *Am. Math. Monthly*, 88, 47–49